%revised 10/4/02 by Matti Lassas
\input amstex
\documentstyle{amsppt}
\magnification=1000
\NoBlackBoxes
%\define\newtext{}

\def\R{\Bbb R}

\define\supp{\text{supp }}

\define\dist{\text{dist }}

\redefine\exp{{\text{\rm exp}}}

\def\fp{\flushpar}

\def\g{\gamma}

\def\L{\Lambda}

\def\O{\Omega}

\def\rta{\rightarrow}

\def\br{\Bbb R}

\def\l{\lambda}

\def\po{\partial\O}

\def\supp{\text{supp}}
\def\D{\Delta}
\def\Dr{\D_\rho}

\def\hra{\hookrightarrow}

\def\back{\backslash}

\def\dv{\vec{\delta}}

\def\path{\partial_\theta}
\def\opt{\ltr}

\def\CD{\Cal C\Cal D}
\def\xpr{X^{p,r}}
\def\codim{\text{ codim}}
\def\p{\partial}
\def\cl{C^l_*}
\def\clsm{\cl S^{m,\vec\delta}}
\def\cli{\cl I^{m,\vec\delta}(H)}
\def\clk{\cl I^{-k,\vec\delta}(H)}
\def\clkz{\cl I^{-k,\vec\delta_0}(H)}
\def\ltr{\langle\theta\rangle}

\NoRunningHeads
\topmatter
\title
The Calder\'on problem for conormal
potentials, I: Global uniqueness and reconstruction
\endtitle
\author
 Allan Greenleaf, Matti Lassas and Gunther Uhlmann
\endauthor
%\date{June 25, 2002}\enddate
\thanks
This work was supported in part by
 National Science Foundation grants DMS-9877101, DMS-0070488
and Finnish Academy project 172434. The
third author was also supported by a John Simon Guggenheim Fellowship.
The second and third author are grateful for the hospitality of the
Mathematical Sciences Research Institute in Berkeley, CA.
\endthanks

\address
Department of Mathematics, University of Rochester, Rochester, NY 14627
\endaddress
\address
Rolf Nevanlinna Institute, P.O. Box 4, 00014
University of Helsinki, Finland
\endaddress
\address
Department of Mathematics, University of Washington, Seattle, WA 98195
\endaddress

\endtopmatter
\document
\head
{\bf 1. Introduction}
\endhead
The goal of this paper is to establish global uniqueness and obtain
reconstruction, in dimensions
$n\ge3$, for the Calder\'on problem in the class of potentials conormal to
a
smooth submanifold $H$ in $\R^n$. In the case of hypersurfaces, the
potentials
considered here may have any singularity weaker than that of the delta
function
$\delta_H$ on the
hypersurface $H$; in general, these potentials correspond to
conductivities
which are
in $C^{1+\epsilon}$ and thus fail to be covered by previously known
results.

Let $\O\subset\R^n$ be a bounded Lipschitz domain, $H\subset\O$ a smooth
submanifold of codimension $k$, and $q\in I^{\mu}(H)$ a real conormal
distribution of
order
$\mu$ with $\mu<1-k$. Thus, if $H=\{x: F_j(x)=0, 1\le j\le k\}$ is a local
representation of
$H$ by means of  defining functions with $\{\nabla F_j: 1\le j\le k\}$
linearly
independent on
$H$, then locally
$q(x)$ has the Fourier integral representation
$$q(x)=\int_{\R^k} e^{i\sum_j F_j(x)\cdot\theta_j} a(x,\theta)
 d\theta,\quad a\in
S^{\mu}_{1,0},\tag1.1$$
where $S^{\mu}_{1,0}$ denotes the standard class of symbols of order $\mu$
and
type $(1,0)$ on $\R^n\times (\R^k\back 0)$. (Here, we use the order
convention
of [12] rather than [16].) A general element
$q\in I^{\mu}(H)$ is a locally finite sum of such expressions. We assume
throughout that
$\supp(q)$ is compact in
$\O$. If $-k<\mu<0$,
then $q$ satisfies $|q(x)|\le C\cdot\dist(x,H)^{-k-\mu}$, so that $q\in
L^{\frac{k}{k+\mu}-\epsilon}(\O),\forall\epsilon>0$, and no better in
general;
in particular, a general element of $I^\mu(H)$ is unbounded. For
comparison, surface measure $\delta_H\in I^0(H)$ and, in the hypersurface
case, a
Heaviside discontinuity across
$H$ belongs to $I^{-1}(H)$.

Rather than working with the Dirichlet-to-Neumann map, $\L_q$, we state
our
main results
in terms of the Cauchy data, $\CD_q$, of sufficiently regular solutions of
the
Schr\"odinger equation
$$(\D+q(x))u(x)=0\hbox{ on }\O.\tag1.2$$
This is more flexible, since $\CD_q$ can be defined for potentials $q$ for
which
$\L_q$ is either not defined (for example, if $\lambda=0$ is a Dirichlet
eigenvalue) or is not known to be defined (due to the low regularity of
$q(x)$); it is perhaps more natural as well.
\medskip
It will be convenient to write $\mu=\nu-k$. Assume that
$\nu_0(k)<\nu<1$, where
$\nu_0(k)\buildrel{def}\over{=}max(\frac23,1-\frac{k}4)$. Fix $p$ and $r$
satisfying
$$2\le
r<\frac{k}{2(1-\nu)}\buildrel{def}\over{=}
r_0(k,\nu)<p_0(k,\nu)\buildrel{def}\over{=}\frac{2k}{k-\nu}<p<\infty.\tag1.3$$
(If
$\nu\le\nu_0(k)$, just pick $p$ and $r$ for some $\nu'>\nu_0(k)$.) Fixing
a
smooth function $\psi\in C^\infty(\R^n), \psi\equiv 1$ near
$\partial\O$ and $int\{\psi=0\}\cap\O\ne \phi$, define the norm
$$||f||_{\xpr} =||f||_{L^p(\O)}+||\D
f||_{L^{p'}(\O)}+||\psi f||_{W^{2,r}(\O)},\tag1.4$$
where $p'$ is the dual exponent to $p$ and $W^{2,r}$ is the standard
Sobolev
space
of $f\in\Cal D'(\O)$ having two derivatives in $L^r(\O)$. Set
$$\xpr(\O)=\{f\in\Cal D'(\O): ||f||_{\xpr} <\infty\}\tag1.5$$
and note that the Schr\"odinger operator $\D+q$ maps $\xpr(\O)\rta
L^{p'}(\O)$
continuously. We denote throughout this paper by $n$
the unit outer normal to $\O$.
\medskip
\proclaim{Definition} For a potential $q\in I^\mu(H)$ with
$H\cap\supp(q)\subset int\{\psi=0\}$, the Cauchy data of the
Schr\"odinger operator $\D+q$ relative to $\xpr(\O)$ is
$$\CD_q=\CD_q^{p,r}=
\Big\{(u|_{\po},\frac{\partial u}{\partial n}|_{\po}): u\in\xpr(\O),
(\D+q)u=0\hbox{ on }\O\Bigr\}.\tag1.6$$
\endproclaim
\medskip
By Sobolev embedding, $\CD_q$ is a subspace of $W^{2-\frac1r,r}(\po)\times
W^{1-\frac1r,r}(\po)$.
Observe that if the Dirichlet-to-Neumann map $\L_q$ {\it is} defined (on
$W^{2-\frac1r,r}(\po)$, say), then
$\CD_q$
is simply its graph. We will construct certain
nontrivial exponentially growing solutions $u\in
\xpr(\O)$, so that, for the potentials considered, $\CD_q$ is in fact
nontrivial. We can now state our first result.

\bigskip
\proclaim{Theorem 1} Suppose that for $j=1,2$, $H_j\subset\subset\O$ are
submanifolds of codimension $k_j$. Suppose further that $q_j\in
I^{\mu_j}(H_j)$
are real potentials  with $\nu_0(k_j)-k_j<\mu_j<1-k_j$ and
$\supp(q_j)\subset\subset\O$. Let $p,r$ satisfy $2\le
r<min(r_0(k_1,\nu_1),r_0(k_2,\nu_2))$ and
$max(p_0(k_1,\nu_1),p_0(k_2,\nu_2))<p<\infty$ and suppose that $\psi\equiv
0$ on
a neighborhood
of $H_1\cup H_2$.
Then $\CD_{q_1}=\CD_{q_2}$ relative to $\xpr(\O)$ implies that
$q_1=q_2$ on
$\O$.
\endproclaim
\bigskip

We also show that under the same assumptions as in Theorem 1 for
the potential
we have a reconstruction procedure, that is we can reconstruct $q$
from $\CD_q$ (see Theorem 2 in section 3 for more details).

Global uniqueness was established
in [32] for $n\ge 3$ (for smooth potentials) and [26] for $n=2$; for $n\ge
3$
this was extended to $q\in L^\infty$ in [27]. The regularity  was
further lowered to $q\in L^{\frac{n}2}$ in unpublished work of R. Lavine
and A. Nachman and  to potentials of small norm in the Fefferman-Phong
class in [4].
Note that for
$-\frac{n-2}n k\le\mu<0, k<\frac{n}2$, a general element of $I^\mu(H)$
fails to
be in
$L^{\frac{n}2}(\O)$.

\medskip
The isotropic conductivity problem, where one considers the
Dirichlet-to-Neumann map for $L_\g=\nabla(\g\cdot\nabla)$, can be
reduced to the Schr\"odinger problem via the substitution
$q=-\frac{\D(\g^{\frac12})}{\g^\frac12}$, and thus
the analogue
of the Theorem holds
for conductivities $\g_j\in I^{-k-1-\epsilon}(H_j)\hra
C^{1+\epsilon}(\overline\O),\forall 0<\epsilon<1$. Currently, the
best general global uniqueness result known for $n\ge 3$ is for
$\g\in C^{\frac32}$,  proved in [28], building on [2] and using the
general argument of [32], while the best known result for $n=2$ is
$\g\in W^{1,p}(\O)$, $p>2$, proved in [3] using the
$\overline\partial$ technique of [1,25,26]. Global
uniqueness for piecewise-analytic conductivities was
proven [20], and special types of  jump
discontinuities  were treated in  [17].

\bigskip
Here, we will follow the general argument of [32], although employing a
different
integral identity so as to avoid difficulties when applying Green's
Theorem.
It is this identity that makes $\xpr$ a convenient space for the problem;
indeed, both sides of $\int_\O \D u\cdot v- u\cdot\D v dx=\int_{\po}
\partial_n u\cdot v-u\cdot\partial_n v d\sigma$ are continuous with
respect
to
$||\cdot ||_{\xpr}$ and thus Green's Theorem holds for $u,v\in\xpr(\O)$.
\bigskip
We now start the proof of Thm. 1, so as to motivate the technicalities
that
follow.

 Given a submanifold $H$ of codimension $k$ and a potential
$q\in I^{\mu}(H)$ with $\mu<1-k$, we will construct exponentially growing
solutions
of (1.2) belonging to
$\xpr(\O)$, of the form $v(x)= e^{\rho\cdot x}(1+\psi(x,\rho))$,  with
$\rho\in\Bbb C^n$ satisfying $\rho\cdot\rho=0$.
Let $v_1\in \xpr(\O)$ be a solution of
$(\D+q_1)v_1=0$. By the hypothesis of Thm. 1, there is a solution
$v_2\in\xpr(\O)$
of $(\D+q_2)v_2=0$ with
$$v_2|_{\po}=v_1|_{\po}\quad\hbox{ and }\quad \frac{\partial v_2}{\partial
n}|_{\po}=
\frac{\partial v_1}{\partial n}|_{\po}.\tag1.7$$
Let $w_2\in\xpr(\O)$ be any other solution to  $(\D+q_2)w_2=0$. Then,
$$(\D+q_2)(v_1-v_2)=(\D+q_2)v_1=(\D+q_1+(q_2-q_1))v_1=(q_2-q_1)v_1,$$
so that
$$\eqalign{\int_\O(q_2-q_1)v_1w_2 dx=&\int_\O (\D+q_2)(v_1-v_2)\cdot w_2
dx\cr
=&\int_\O (\D+q_2)(v_1-v_2)\cdot w_2-(v_1-v_2)\cdot (\D+q_2)w_2 dx\cr
=&\int_\O \D(v_1-v_2)\cdot w_2 - (v_1-v_2)\cdot\D w_2 dx\cr
=& \int_{\po} \frac{\partial}{\partial n}(v_1-v_2)\cdot w_2 -
(v_1-v_2)\cdot
\frac{\partial}{\partial n} w_2 d\sigma\cr
=& 0,}\tag1.8$$
where the application of Green's Theorem is valid since $v_1-v_2$ and $
w_2\in
\xpr$ and the last equality holds by (1.7). If we carry this out for the
solutions
$v_1$ and
$w_2$ constructed below for complex frequencies $\rho_1$ and $\rho_2$
satisfying
$\rho_1+\rho_2=-i\xi$, with $\xi\in\R^n\back 0$, then we have, as in [32],

$$\eqalign{0=& \int_\O (q_1-q_2) e^{(\rho_1+\rho_2)\cdot x}
\Bigl(1+\psi_1(x,\rho_1)\Bigr)\Bigl(1+\psi_2(x,\rho_2)\Bigr) dx\cr
=&\widehat{(q_1-q_2)}(\xi) +\int_\O e^{-i\xi\cdot x}
(q_1-q_2)(\psi_1+\psi_2+\psi_1\psi_2)dx}.\tag1.9$$
If one can do this for pairs
$(\rho_1,\rho_2)$ with $|\rho_j|\rta\infty$ and show that the last
integral $\rta 0$ as $|\rho|\rta\infty$, then
$\hat{q_1}(\xi)=\hat{q_2}(\xi)$; doing this for all
$\xi\in\R^n$ will finish the proof of Thm. 1.

As is well known, $v(x)= e^{\rho\cdot x} (1+\psi(x))$ is a solution
of the Schr\"odinger equation iff $\psi$ is a solution of
$$(\D_\rho+q)\psi=-q(x)\hbox{ where }
\D_\rho=\D+2\rho\cdot\nabla.\tag1.10$$
We will show in Prop. 2.6 that (1.10) is uniquely solvable, with some
decay in
$|\rho|$, in a
Banach space of finite-regularity conormal distributions associated with
$H$, yielding
exponentially growing solutions $v_j\in\xpr$ to (1.2) which allow the
argument above to
be carried out. In \S3, this result is extended to a hybrid global space,
and
this is applied to obtain reconstruction of the potential from the
Cauchy data, following the general argument of [25]. Finally, in \S4 we
show that uniqueness can fail in a weak formulation of the problem for
potentials with very strong singularities on a hypersurface,
with  blow-up rates corresponding
  to those of distributions conormal of order greater than 1 for $H$.

We would like to thank Steve McDowall for valuable discussions, and the
referee
for pointing out an error in the original version of this paper.
\bigskip
\head
{\bf 2. Uniqueness for conormal potentials}
\endhead

As described in the Introduction, to prove Thm.1, it suffices to construct
exponentially growing solutions to (1.2) of the form $v(x)=e^{\rho\cdot
x}(1+\psi(x,\rho))$ for $\rho\cdot\rho=0,|\rho|\rta\infty$ so that the
second
integral in (1.9) tends to 0 as $|\rho|\rta\infty$. To do this for
potentials
$q\in I^\mu(H)$, the standard space of (infinite-regularity) conormal
distributions of order
$\mu$ associated with the codimension $k$ submanifold
$H$, we will also need to formulate Banach spaces of finite-regularity
conormal
distributions in
$\R^n$ associated with $H$.  Rather than working in unnecessary
generality, we
will restrict ourselves to the spaces  needed here; unlike [23],[22],
where
several other types of finite-regularity conormal spaces are defined using
iterated regularity with respect to Lie algebras of tangent
vector fields, we impose the finite-regularity assumption directly on the
symbols in the oscillatory representations of the distributions, using
symbol
classes modelled on those of [33].

For $l\in\R$, let $\cl$ denote the Zygmund space of order $l$ on $\R^n$
[33].
Thus, if $\psi_0(D)+\sum_{i=1}^\infty \psi_i(D)=I$ is a Littlewood-Paley
decomposition, with $\psi_i(\xi)=\psi_1(\frac\xi{2^{i-1}}), i\ge 1$, then
$$||u||_{\cl} =\sup_i 2^{li} ||\psi_i(D)u||_{L^\infty(\R^n)}.$$
Recall that if $l\ge 0, l\notin\Bbb Z$, then $\cl= C^{[l],l-[l]}(\R^n)$.

Now
fix an order
$m\in\R$, an $N\in \Bbb N$, and a sequence
$\dv=(\delta_1,\delta_2,...,\delta_N)$ of numbers
$0\le \delta_j\le 1, 1\le j\le N$. For any multi-index $\alpha\in \Bbb
Z_+^k$,
let $\delta(\alpha)=\sum_{j=1}^{|\alpha|} \delta_j$, setting $\delta(0)=0$
for
convenience, and $|\dv|=\sum_{j=1}^N \delta_j$.
\bigskip
\proclaim{Definition 2.1} (i)
$$\clsm (\R^n\times\R^k)=\Bigl\{a(x,\theta):
|| \path^\alpha a(\cdot,\theta)||_{\cl}\le
C_\alpha(1+|\theta|)^{m-\delta(\alpha)}, \forall  |\alpha|\le
N\Bigr\}$$ and
$$||a||_{\clsm}=\max_{0\le\alpha\le
N} (1+|\theta|)^{
-m+\delta(\alpha)} ||
\path^\alpha a(\cdot,\theta)||_{\cl}.$$
(ii) If $H$ is a smooth codimension $k$ submanifold with compact closure,
then
$\cli$ is the space of locally finite sums of distributions of the form
$u(x)=\int_{\br^k} e^{iF(x)\cdot\theta}
a(x,\theta)
d\theta$
with $a\in\clsm$, where $F(x)=(F_1(x),\dots,F_k(x))$ are local
defining functions for $H$.
\endproclaim
\fp{\bf Remarks}. (1) $\clsm$ is a Banach space with respect to the norm
defined in
(i), and $\cli$ inherits this structure.
\fp (2) $\path: \clsm\rta \cl S^{m-\delta_1,\vec\delta'}$, with
$\vec\delta'= (\delta_2,\delta_3,...,\delta_N)$ and
$\partial_x: \clsm\rta
C_*^{l-1} S^{m,\dv}$ continuously.
\fp (3) $\cli$ is well-defined, since changing the defining functions of
$H$
locally corresponds to a change of variable in $x$, which leaves the
symbol
class invariant.
\fp (4) The usual (infinite-regularity) conormal space $I^m(H)$ has a
continuous
inclusion
with respect to its Fr\'echet structure: $ I^m(H)\hra \cli$ for any $l$,
$N$
and $\dv$.

\bigskip
\proclaim{Proposition 2.2} Let $H\subset \Omega\subset\subset \R^n$, with
$\codim(H)=k$, $l\notin\Bbb Z$ and $\dv=(\delta_1,\dots,\delta_N)$.

 (i) If $l>1, -k\le m<-\frac{k}2, N\ge1, \delta_1>0$ and
$m-|\dv|<-k$, then
$$\cli\hra L^p(\O)\hbox{\text{ continuously for
 all} }
1\le p <\frac{k}{m+k}.\tag2.1$$
\fp (ii) If $l>2$, $N\ge2$, $\delta_1>1-\frac{k}2$ and
$|\dv|>1$, then for each smooth function $\tilde F$ vanishing on
$H$,
$$v\in \clk\implies \tilde F\cdot v\in
W^{1,p}(\Omega)\hbox{
 for }1\le p<\frac{k}{1-\delta_1}\tag2.2$$
and thus, for $G\subset\subset\overline\O\back H$,
$\clk\hra W^{1,p}(G)$ continuously.
\fp (iii) Suppose $l> 3$. For $\max(\frac23,
1-\frac{k}4)\buildrel{def}\over{=}
\nu_0(k)<\nu<1$, set
$$\dv_0=(\nu,1-\nu,2\nu-1,1-\nu,1-\nu,3\nu-2).$$
Then, for any $\tilde F$
vanishing on
$H$,
$$v\in \clkz\implies \tilde F^3\cdot v\in
W^{2,r}(\Omega)\hbox{
 for }1\le r<\frac{k}{2(1-\nu)}\tag2.3$$
and thus, for $G\subset\subset\overline\O\back H$,
$\clkz\hra W^{2,r}(G)$ continuously.
\endproclaim
\fp Proof. (i) Since $L^p$ and $\cli$ are diffeomorphism-invariant, it
suffices to
assume that, with respect to coordinates
$x=(x',x'')\in\R^{n-k}\times\R^k$,
$$H=\{x''=0\}\quad\hbox{\text and }\quad u(x)=\int_{\R^k}
e^{ix''\cdot\theta}
a(x,\theta)
d\theta,\quad a\in \clsm.$$

We then have
$$u(x',x'')=\int
e^{ix''\cdot\theta}[a(x',0,\theta)+\sum_{0<|\alpha''|<[l]}
\frac1{(\alpha'')!}\partial^{\alpha''}_{x''}a(x',0,\theta)(x'')^{\alpha''}
+R_{[l]}(x',x'',\theta)]
d\theta,$$
where
$$R_{[l]}(x,\theta)=\sum_{|\alpha''|=[l]} b_{\alpha''}(x,\theta)
(x'')^{\alpha''} \hbox{ with } |b_{\alpha''}(x,\theta)|\le
C(1+|\theta|)^{m}, \forall \alpha''.$$
Since $a(x',0,\cdot)\in
L^q(\R^k_\theta),\forall q>-\frac{k}{m}$ and
$-\frac{k}{m}<2$, the
Hausdorff-Young inequality implies that $\int e^{ix''\cdot\theta}
a(x',0,\theta)
d\theta\in L^{q'}(\R^k_{
x''}),\forall 2\le q'<\frac{k}{m+k}$
, uniformly in $x'$, and thus
belongs to
$L^p(\O)$ for all $2\le p<\frac{k}{m+k}$; since it is compactly supported,
the range
is in fact $1\le p<\frac{k}{m+k}$. For the second term, note that for each
function
$a_{\alpha''}(x',\theta)\buildrel{def}\over{=}\frac1{(\alpha'')!}\partial^{\alpha''}a(x',0,\theta)$,
$$\int e^{ix''\cdot\theta} a_{\alpha''}(x',\theta)( x'')^{\alpha''}
 d\theta =
\int ( \frac1{i}\partial_\theta)^{\alpha''}(e^{ix''\cdot\theta})
a_{\alpha''}(x,\theta) d\theta =
$$
$$=
\int e^{i x''\cdot\theta}
(\frac{-1}{i}\partial_\theta)^{\alpha''} (a_{\alpha''}(x,\theta))
d\theta,$$
whose amplitude is $\le C(1+|\theta|)^{m-\delta(\alpha'')}$, and
hence is treated by the same argument as for $a_0=a(x',0,\theta)$. In the
final term, we integrate by parts:
$$\int e^{ix''\cdot\theta}  b_{\alpha''}(x,\theta) (x'')^{\alpha''}
d\theta  =
\int e^{i x''\cdot\theta}
(i\partial_{\theta})^{\alpha''} b_{\alpha''}(x,\theta) d\theta,$$
and since $|\partial_{\theta}^{\alpha''} b(x,\theta)|\le
C(1+|\theta|)^{m-|\dv|}\in L^1(\R^k_\theta)$, uniformly in
$x'$, this yields a bounded function of $x\in\O$.
\medskip
\fp (ii) If $v\in \clk$, then $v\in L^p,\forall
p<\infty$,
by part (i). From
$$v(x)=\int e^{ix''\cdot\theta} a(x,\theta) d\theta,\quad
a\in \cl S^{-k,\dv},$$
one finds
$$\eqalign{\partial_{x_j} v(x)=&\int
e^{ix''\cdot\theta}
(i\theta_j a+\partial_{x_j}a) d\theta\cr
\in& \cl I^{1-k,\dv}(H) + C_*^{l-1} I^{-k,\dv}(H)}$$
for $n-k+1\le j\le n$; if $1\le j\le n-k$, only the second term is
present.
The second term is covered by part (i) and hence
is in $L^p,\forall p<\infty$. If
the first term is multiplied by some $x_{j_0}, n-k+1\le j_0\le n$  and
integrated by parts, it becomes an element of
$\cl I^{1-k-\delta_1,\dv'}(H)$, with $\dv'=(\delta_2,\dots,\delta_N)$,
which by
(2.1) is in
$L^p, 1\le p<\frac{k}{1-\delta_1}$. Since any $\tilde F$ vanishing on $H$
can be
represented as a linear combination of $x_{j_0}$'s with smooth
coefficients,
$\tilde F(x)\nabla v\in L^p$ and so
$v\in W^{1,p}(G)$ for any set
$G$ on which
$|F|$ is bounded below.
\medskip
\fp (iii) By  (i) and (ii) above, both $v$ and $\tilde F\cdot\nabla v\in
L^r$.
Now, arguing as in (ii), for $n-k+1\le j,j'\le n$,
$$\eqalign{\partial^2_{x_jx_{j'}}v=&\int e^{i
x''\cdot\theta}\Bigl(-\theta_j\theta_{j'}a+i(\theta_j
a_{x_{j'}}+\theta_{{j'}}a_{x_j})+a_{\theta_j \theta_{j'}}\Bigr) d\theta\cr
\in&\quad\cl I^{2-k,\dv_0}(H)+ C_*^{l-1}I^{1-k,\dv_0}(H)
+ C_*^{l-2}I^{-k,\dv_0}(H),}$$
with simpler expressions if one or both of $j$ or $j'$ is $\le n-k$.
By (i), the last term is in $L^r,\forall  r<\infty$, if $l>3,N\ge
1$. Integrating by parts and using (ii), $x_{j_0}$ times the second term
is in
$L^r, 1\le r<\frac1{1-\delta_1}=\frac1{1-\nu}$ if $l> 2,N \ge 2$.
As
for the first term, $x_{j_0}x_{j_1}x_{j_2}$ times it is seen, after
integrating
by parts three
times, to be in
$\cl I^{2-k-(\delta_1+\delta_2+\delta_3),\dv_0''}(H)
=\cl I^{2-k-2\nu,\dv_0''}(H)$ for $\dv_0''= (1-\nu,1-\nu,3\nu-2)$,
which, by (i), $\hra L^r(\O),\forall 1\le r<\frac{k}{2(1-\nu)}$, if $l>
2$, since $-k\le 2-k-2\nu<-\frac{k}2$ if $1-\frac{k}4<\nu\le 1$
and $2-k-|\dv|=2-3\nu-k<k$ if $\nu>\frac23$. Thus, for any $\tilde F$
vanishing on $H$,
$\tilde F^3\cdot\clkz\hra W^{2,r}(\O)$ if $1\le
r<\frac{k}{2(1-\nu)},l> 3$.
\qed
\bigskip
We also have
\bigskip
\proclaim{Proposition 2.3} If $A\in\Psi^r(\R^n)$ is properly supported
with
$r<0$, then for any
$l,m\in\R$ and any $\dv$, and any $\epsilon>0$,
$$A:\cli\rta \cl I^{m+r+\epsilon,\dv}(H)$$
and $||A||$ is bounded by a finite number of semi-norms of $\sigma(A)$ in
$S^r_{1,0}(\R^n\times(\R^n\back 0))$.
\endproclaim
\bigskip
\fp{\bf Proof.}
Write $$Af(x)=\int e^{i(x-y)\cdot\xi} a(x,y,\xi) f(y) d\xi
dy,\quad
a\in S^r_{1,0}$$
and
$$u(y)=\int e^{iF(y)\cdot\theta} b(y,\theta) d\theta,\quad
 b\in\cl S^{m,\dv}.$$
Then
$$\eqalign{Au(x)=&\int e^{i[(x-y)\cdot\xi+F(y)\cdot
\theta]} a(x,y,\xi) b(y,\theta)
d\theta d\xi dy\cr
=&\int e^{iF(x)\cdot\theta}c(x,\theta) d\theta,}$$
where
$$
c(x,\theta)=
\int e^{i[(x-y)\cdot\xi+ (F(y)-F(x))\cdot\theta]}
a(x,y,\xi)
b(y,\theta) d\xi dy\tag2.4
$$
and we need to show that $c=c[b]\in\cl
S^{m+r+\epsilon,\dv},\forall\epsilon>0$.
 Furthermore, since $H$ is compact, we
may assume that the symbol $a(x,y,\xi)$ vanishes for $x,y$ outside some
compact set. 

Let $\ltr=(1+|\theta|^2)^{\frac12}$. Setting
$b_i(y,\theta)=\psi_i(D_y)(b(y,\theta))$ and
$c_j(x,\theta)=\psi_j(D_x)(c(x,\theta))$, we have
$$||\path^\alpha
b_i(\cdot,\theta)||_{L^\infty}\le C
2^{-li}\ltr^{m-\delta(\alpha)},\forall |\alpha|\le N,$$ and want to
prove that
$$||\path^\alpha
c_j(\cdot,\theta)||_{L^\infty}\le C_\epsilon
2^{-lj}\ltr^{m+r+\epsilon-\delta(\alpha)},\forall |\alpha|\le N.$$

We first
consider the case of $\alpha=0$.  Let $\{\tilde\psi_i\}_{i=0}^\infty,
\{\tilde{\tilde\psi}_i\}_{i=0}^\infty$ be
 bounded families in $S^0_{1,0}$ with $\tilde\psi_i=1$ on $\supp(\psi_i)$
and $\tilde{\tilde\psi}_i=1$ on $\supp(\tilde\psi_i)$; then, one can write
$c_j=T_{ij}(b_i)$, where the operator $T_{ij}$ has the Schwartz kernel
$$K_{ij}(x,y;\theta)=\int
e^{i[(x-z)\cdot\zeta+(z-w)\cdot\xi+(F(w)-F(z))\cdot\theta
+(w-y)\cdot\eta]} \psi_j(\zeta) a(z,w,\xi)\tilde\psi_i(\eta) d\zeta dz
d\xi dw
d\eta.$$
It suffices to show that $\sum_{i=0}^\infty ||T_{ij}(b_i)||_{L^\infty} \le
C_\epsilon 2^{-lj}\ltr^{m+r+\epsilon}$.

As in the proof of Prop. 2.2, we may assume that $F(x)=x''$, where
$x=(x',x'')\in\R^{n-k}\times\R^k$; for $\theta\in\R^k$, let
$\theta^*=(0,\theta)=DF^*(x)(\theta),\forall x$. Thus, the phase function
of $K_{ij}(x,y;\theta)$ is
$$\phi=(x-z)\cdot\zeta+(z-w)\cdot\xi +(w''-z'')\cdot\theta
+(w-y)\cdot\eta.$$
Let $\chi\in C_0^\infty(\R), \chi(t)=1$ for $|t|\le
1/2,\supp(\chi)\subset\{|t|\le 3/4\}$. Write
$$\eqalign{K_{ij}(x,y;\theta)=&\int
e^{i\phi}\psi_j(\zeta)\tilde{\psi}_i(\eta)[\chi
(\frac{|\xi-\zeta-\theta^*|}{|\theta|})+
(1-\chi)(\frac{|\xi-\zeta-\theta^*|}{|\theta|})]\times\cr
&\quad\quad 
a(z,w,\xi) d\zeta d\eta d\xi dz
dw\cr
=& K_{ij}^0(x,y;\theta) + K_{ij}^\infty(x,y;\theta)}
$$
with $T_{ij}=T_{ij}^0+ T_{ij}^\infty$ the corresponding operator
decomposition.

$K_{ij}^\infty$ may be analyzed by noting that $|\xi-\zeta-\theta^*|\ge
c\max(|\xi-\zeta|,|\theta|)$ on $\supp(1-\chi)$, and thus, using
$$(\frac{\xi-\zeta-\theta^*}{i|\xi-\zeta-\theta^*|^2}
\cdot\nabla_z)(e^{i\phi})=e^{i\phi},$$ 
we may integrate by parts $n+M$ times in $z$ and then integrate in $\xi$ to
obtain
$$K_{ij}^\infty(x,y;\theta)=\int e^{i[(x-z)\cdot\zeta
+(w''-z'')\cdot\theta +(w-y)\cdot\eta]} \psi_j(\zeta)\tilde{\psi}_i(\eta)
A(z,w;\theta) d\zeta d\eta dz dw$$
with $|\partial_z^\alpha \partial_w^\beta A(z,w;\theta)|\le
C_{\alpha,\beta,M}\ltr^{-M},\forall \alpha,\beta\in\Bbb Z_+^n, M\in\Bbb
N$. If $\max(2^i,2^j)<\frac14\ltr$,  we then simply integrate in all
variables to obtain $|K_{ij}^\infty|\le c2^{n(i+j)}\ltr^{-M}$. If
$\ltr<\frac14\min(2^i,2^j)$, then $|\zeta+\theta^*|\ge c|\zeta|\ge c2^j$
and $|\eta+\theta^*|\ge c|\eta|\ge c2^i$, so we may first integrate by
parts in $z$ and $w$ to obtain $|K_{ij}^\infty|\le
c2^{-M'(i+j)}\ltr^{-M}$. 

In the following,  we denote $a\sim b$ if either
$a,b<4$ or $\frac 14 a\leq b\leq 4a$.

If $\ltr\sim 2^i>2^{j+1}$ or $\ltr\sim 2^j>
2^{i+1}$, we may integrate by parts in $z$ or $w$ alone to obtain
$|K_{ij}^\infty|\le c 2^{-M'i}\ltr^{-M}$ or $|K_{ij}^\infty|\le c
2^{-M'j}\ltr^{-M}$, respectively, $\forall M'$. Finally, if $\ltr\sim
2^i\sim 2^j$, we simply integrate and get $|K_{ij}^\infty|\le c
2^{n(i+j)}\ltr^{-M}\le c\ltr^{-(M-2n)}$. Using $||b_i||_{L^\infty}\le c
2^{-li}\ltr^m$ and summing in $i$ yields $\sum_{i=0}^\infty
||T_{ij}^\infty(b_i)||_{L^\infty}\le c 2^{-M'j}\ltr^{-M},\forall M,M'$.

On the other hand, $K_{ij}^0$ may be analyzed using stationary phase in
$\xi$ and $w$: First rewrite 
$$\eqalign{&K_{ij}^0(x,y;\theta)=\int
e^{i[(x-z)\cdot\zeta+(z-y)\cdot\eta]}
\psi_j(\zeta)\tilde{\psi}_i(\eta)\times  \cr
&\Bigl[ \int
e^{i|\theta|[w''\cdot\frac{\theta}{|\theta|}-w\cdot
\xi)]}
|\theta|^n a(z,z+w,|\theta|(\xi+\frac \eta{|\theta|}))\chi(|\xi-
\frac{\theta^*}{|\theta|}-\sigma|) d\xi
dw\Bigr]_{\sigma=\frac {\zeta-\eta}{|\theta|}} dzd\zeta d\eta.
}\tag2.5
$$ 

The domain of integration for the inner integral
contains a critical point only if $|\sigma|<2$. In this 
case the unique nondegenerate critical point is
$\xi=\frac{\theta^*}{|\theta|},\ w=0$.
Applying the method stationary phase 
and analysing error terms using [16,Thm.7.7.7]
for $|\sigma|<2$, or just using integration by parts for $
|\sigma|\geq 2$, we 
see that the inner
integral is 
$a(z,z,\eta+\theta^*)\chi(\frac {|\eta-\zeta|}{|\theta|})
 + e(z,\frac{\eta}{|\theta|},\frac{\zeta}{|\theta|},
\eta+\theta^*)$, with $ e\in S^{r-1}_{1,0}((\R^n\times \R^n\times \R^n)\times
\R^n
\back 0)$ 
and
thus
$$
\eqalign{
K_{ij}^0(x,y;\theta)=\int& e^{i[(x-z)\cdot\zeta +(z-y)\cdot\eta]}\psi_j(\zeta)
\tilde\psi_i(\eta)\times \cr
& [a(z,z,\eta+\theta^*)\chi(\frac{|\eta-\zeta|}{|\theta|})+
e(z,\frac{\eta}{|\theta|},\frac{\zeta}{|\theta|},\eta+\theta^*)] 
d\zeta dz d\eta.}
\tag2.6$$
Note that $S^{r-1}_{1,0}$ seminorms of $e$ are uniformly
bounded over compact subsets of 
$\R^n\times \R^n\times \R^n$.
We observe for the phase that $d_z
((x-z)\cdot\zeta +(z-y)\cdot\eta)=-\zeta+\eta$ and
$(d_\zeta+d_\eta)((x-z)\cdot\zeta +(z-y)\cdot\eta)=x-y$ and we may use
these to
integrate by parts, $M$ and $M'$ times, respectively, to obtain, for $i\le
j-2$,
$$\eqalign{|K^0_{ij}(x,y;\theta)|\le&C\int
2^{-Mj}(1+2^i|x-y|)^{-M'}\tilde\psi_j(\zeta)\ltr^r
|\eta|^{-M'}\tilde{\tilde\psi_i}(\eta) d\zeta dz d\eta\cr
\le&C 2^{-(M-n)j} (1+2^i|x-y|)^{-M'} 2^{-(M'-n)i}\ltr^r.}$$
Here, we have used $|-\zeta+\eta|\ge\frac12|\zeta|\ge 2^{j-1}$ and
$$|(d_\zeta+ d_\eta)^{M'}(\psi_j(\zeta)
a(z,z,\eta+\theta^*)\tilde\psi_i(\eta))\le
\cases C|\eta|^{r-M'}\tilde\psi_j(\zeta)\tilde{\tilde\psi_i}(\eta), |\theta|\le c
2^i\cr
C\ltr^r |\eta|^{-M'}\tilde\psi_j(\zeta)\tilde{\tilde\psi_i}(\eta),
|\theta|\ge
c2^i.\endcases$$
Thus, $\int|K^0_{ij}(x,y)|dy \le C 2^{-(M-n)j} 2^{-M'i}\ltr^r$, so that
$$\sum_{i=0}^{j-2} ||T^0_{ij}(b_i)||_{L^\infty}\le C\sum_{i=0}^{j-2}
2^{-(M-n)j}
2^{-M'i}\ltr^r 2^{-li} \ltr^m\le C c^{-lj}\ltr^{m+r}$$
for $M>n+l$. Similarly, for $i\ge j+2$, we obtain the estimate
$$|K^0_{ij}(x,y)|\le C 2^{-(M-n)i} (1+ 2^j|x-y|)^{-M'}
2^{-(M'-n)j}\ltr^r,$$
leading to
$$\sum_{i=j+2}^\infty ||T^0_{ij}(b_i)||_{L^\infty}\le C \sum_{i=j+2}^\infty
2^{-M'j} 2^{-(M-n)i}\ltr^r 2^{-li}\ltr^m\le C
2^{-lj}\ltr^{m+r}$$
for $M'>l, M>n-l$.

To analyze the contributions from $|i-j|\le 1$, we perform in (2.6) an
additional stationary phase in
$z$ and $\zeta$, with  critical point $z=x$, $\zeta=\eta$. Rewriting the
integral as in formula (2.5) we see that,  
$$\eqalign{K^0_{ij}(x,y)= \int e^{i(x-y)\cdot\eta} \psi_j(\eta)
[(a(x,x,\eta+\theta^*)+&e(x,\frac {\eta}{|\theta|},\frac {\eta}{|\theta|},
\eta+\theta^*))\times\cr
&(1+e_2(x,\frac \eta{|\theta|},\eta+\theta^*))
]
\tilde\psi_i(\eta)d\eta}
$$
with $e_2\in S^{-1}_{1,0}((\R^n\times \R^n)\times \R^n\back 0)$.
Since $\psi_j$ are uniformly bounded
in $S^0_{1,0}$, the 
$S^{-1}_{1,0}$ seminorms of $e_2$ are uniformly bounded over compact subsets
of
$\R^n\times \R^n$. 
Then $$
K^0_{ij}(x,y)= e^{-i\theta^*\cdot x}\Bigl[\int
e^{i(x-y)\cdot\eta}[\psi_j(\eta-\theta^*)
h(x,\eta;\, \theta)\tilde{ \tilde\psi}_i(\eta-\theta^*)]\tilde\psi_i(\eta-\theta^*)
d\eta\Bigr]
e^{i\theta^*\cdot y}
$$
where $h(\cdotp,\cdotp;\, \theta)\in S^0_{1,0}(\R^n \times\R^n\back 0)$,
$$
h(x,\eta;\, \theta)=  
(a(x,x,\eta)+e(x,\frac {\eta-\theta^*}{|\theta|},\frac {\eta-\theta^*}
{|\theta|},
\eta))(1+e_2(x,\frac {\eta-\theta^*}{|\theta|},\eta))
$$
is a symbol depending on the parameter $\theta$, $|\theta|>1$.

Note that multiplication by the exponentials does not affect the
$L^\infty\to
L^\infty$ operator norm, while
$$\psi_j(\eta-\theta^*)h(x,\eta;\, \theta)
\tilde{\tilde\psi}_i(\eta-\theta^*)\in\cases 
S^r_{1,0}(\R^n\times \R^n\back 0),\quad
|\theta|\le c
2^j\cr
\ltr^r\cdot S^0_{1,0}(\R^n\times \R^n\back 0),\quad |\theta|\ge c2^j\endcases$$
uniformly in $\theta$. Acting on functions with Fourier transforms
supported in
$B(0,R)$, pseudodifferential operators of order 0 are bounded on
$L^\infty$
with norm $\le c \log R$ [33]. Hence, for $|i-j|\le 1$,
$$||T_{ij}(b_i)||_{L^\infty} \le\cases C j 2^{rj} 2^{-lj}\ltr^m\le
C_\epsilon 2^{-lj}\ltr^{m+r+\epsilon},\quad |\theta|\le c 2^j\cr
C j 2^{-lj}\ltr^{m+r}\le C_\epsilon
2^{-lj}\ltr^{m+r+\epsilon},\quad
|\theta|\ge c 2^j,\endcases$$
as desired.

To handle the derivatives  $\path^\alpha c$, note first that the
application of
$\partial_\theta$ to the right hand side of (2.4) yields two terms. Since
$\partial_\theta b(x,\theta)\in
\cl S^{m-\delta_1,(\delta_2,...)}$, the oscillatory integral
with
amplitude $a(x,y,\xi)\partial_\theta b(x,\theta)$ has $\cl$ norm $\le
C\ltr^{m+r-\delta_1}$, as desired. On the other hand,
as here $F(x)=x''$, we see that if the
derivative
hits the phase, then the amplitude becomes
$i(F(y)-F(x))\cdot a\cdot b
=i(y''-x'')\cdot a\cdot b$. Writing $
y''-x''=-\partial_{\xi''}[(x-y)\cdot\xi+ (F(y)-F(x))\cdot\theta]$, we may
then integrate by parts in $\xi$.
The resulting amplitude, $i\partial_\xi
a\cdot b$ gives a term whose $\cl$ norm is $\le C\ltr^{m+r-1}$.
Higher derivatives $\path^\alpha c, |\alpha|\le N$, are handled similarly.
\qed
\bigskip

To obtain the boundedness of $M_q$ on some of the finite-regularity
conormal
spaces,
we will need the following; similar results are in [13].
\bigskip
\proclaim{Lemma 2.4} Let $a,b\in C(\R^n\times\R^k)$ satisfy
$$||a(\cdot,\theta)||_{\cl}\le C\ltr^m\quad\hbox{ and }\quad
||b(\cdot,\theta)||_{C_*^{l'}}\le C\ltr^{m'}$$
with $l,l'>0, l,l'\notin\Bbb N$ and $m+m'<-k$.
Then the partial convolution $a*_{\R^k}b$ satisfies
$$||a*_{\R^k}b(\cdot,\theta)||_{C_*^{l''}}\le
C\ltr^{m''},$$
with $l''=\min(l,l')$ and $m''=\max\bigl(
(m+k)_{\tilde+}+m',m+(m'+k)_{\tilde+},m+m'+k)$. Here,
$t_{\tilde+}=t_+=\max(t,0)$ if $t\ne 0$ and $t_{\tilde+}=\epsilon$ for
$t=0$, with
$\epsilon>0$ arbitrarily
 small.
\endproclaim
\bigskip
\fp{\bf Sketch of proof.}
This follows easily by decomposing
$$\eqalign{(a*_{\R^k}b)(x,\theta)=&\int_{\R^k} a(x,\sigma)
b(x,\theta-\sigma)
d\sigma\cr
=&(\int_{|\sigma|\le\frac{|\theta|}2}
+\int_{|\theta-\sigma|\le\frac{|\theta|}2}
+\int_{{|\sigma|\ge\frac{|\theta|}2}
,\ {|\theta-\sigma|\ge\frac{|\theta|}2}}
) a(x,\sigma) b(x,\theta-\sigma)
d\sigma
.}$$
and estimating each of the three terms using  the fact that $\cl\times
C_*^{l'}\hra C_*^{l''}$ continuously.

This is then used in the proof of
\bigskip

\proclaim{Proposition 2.5} Let $H\subset\R^N$ be codimension $k$. For
$q\in
I^\mu(H)$,
let $M_q$ denote multiplication
by $q(x)$. Suppose $\frac23-k\le \mu<1-k$. Then, for any $l>0, l\notin\Bbb
N$,
$$M_q: \clkz\rta
\cl I^{\tilde{\mu},\dv_0}(H)\quad\hbox{continuously}\tag2.7$$
where, if we write $\mu=\nu-k$ with $\frac23\le\nu<1$, $\dv_0$ is as in
Prop.
2.2(iii):
$$\dv_0=(\nu,1-\nu,2\nu-1,1-\nu,1-\nu,3\nu-2)\tag2.8$$
and $\tilde\mu=\mu+\epsilon$ for any $0<\epsilon<-\mu$.
\endproclaim
\bigskip

\fp{\bf Proof.} As before, we can assume that
$H=\{x''=0\}$,
$$u(x)=\int_{\R^k} e^{i x''\cdot\theta} b(x,\theta) d\theta,\quad b\in
\clkz$$
and
$$q(x)=\int_{\R^k} e^{ix''\cdot\theta} a(x,\theta) d\theta,\quad a\in
S^\mu_{1,0}\hra C^{l'} S^{\mu,(1,1,\dots)}$$
for any $l'>0$.

Hence, $M_qu(x)=\int e^{ix''\cdot\theta} (a*b)(x,\theta) d\theta$, where
$*$
denotes
the $k$-dimensional convolution in the $\theta$ variable.
By Lemma 2.4,
$$||a*b(x,\theta)||_{\cl}\le
C\ltr^{\max(\mu,\mu+\epsilon,\mu)}=C\opt^{\mu+\epsilon},\forall
\epsilon>0.$$
Since $\path a\in C_*^{l'}S^{\mu-1,(1,1,\dots)}$ and $\mu-1<-k$,
$$||\path(a*b)(x,\theta)||_{\cl}=||(\path a)*b (x,\theta)||_{\cl}\le
C\opt^{\max(-k,\mu-1+\epsilon,\mu-1)}=C\opt^{-k}$$
(for $\epsilon<-\mu$), which gives a gain of $\ge \mu-(-k)=\nu$,
consistent
with $\delta_1=\nu$. (The additional gain of $\epsilon$ we choose to
ignore.)
Noting that $\path b\in \cl S^{-k-\nu,(1-\nu,2\nu-1,\dots,3\nu-2)}$ by
Remark 2 above, we have
$$||\path^2(a*b)||_{\cl}=||(\path a)*(\path b)||_{\cl}\le
C\opt^{\max(-k-\nu,\nu-k-1,-k-1)}=C\opt^{\nu-k-1}$$
since
$\nu\ge\frac12$, which is a gain of $\delta_2=-k-(\nu-k-1)=1-\nu$. Since
$\path^2a\in C_*^{l'} S^{\mu-2,(1,1,\dots)}$ and $\path b$ is as noted
above,
$$||\path^3(a*b)||_{\cl}=||(\path^2a)*(\path b)||_{\cl}\le
C\opt^{\max(-k-\nu,\nu-k-2)}=C\opt^{-k-\nu},$$ which gives a gain of
$\delta_3=\nu-k-1-(-k-\nu)=2\nu-1$. Since $\path^2a\in C_*^{l'}
S^{\mu-2,(1,1,\dots)}$ and
$\path^2 b\in \cl S^{-k-1,(2\nu-1,1-\nu,...)}$,
$$||\path^4(a*b)||_{\cl}=||(\path^2a)*(\path^2 b)||_{\cl}\le
C\opt^{\max(-k-1,\nu-k-2,2\nu-k-3)}=C\opt^{-k-1},$$
which
gives a gain of $\delta_4=-k-\nu-(-k-1)=1-\nu$.
Continuing in this fashion, we may estimate
$$||\path^5(a*b)||_{\cl}=||(\path^2a)*(\path^3b)||_{\cl}\le
C\ltr^{\nu-k-2},$$
which is consistent with $\delta_5=1-\nu\hbox{ if }\nu\ge \frac23$,
and

$$||\path^6(a*b)||_{\cl}=||(\path^3a)*(\path^3b)||_{\cl}\le
C\ltr^{-2\nu-k},$$
which is consistent with  $\delta_6=3\nu-2$.
 The $x$ derivatives, lowering the Zygmund space index and not involving
any gain in
$\ltr$, are handled in the
obvious way. Hence, $a*b\in \cl S^{\mu+\epsilon,\dv_0}$. \qed
\bigskip

Now recall some facts concerning the
Faddeev Green's function [10], $G_\rho$.  As is well-known (see,e.g.,
[32], where this is used implicitly), the families
$\{|\rho|G_\rho:\rho\cdot\rho=0\}$ and $\{ G_\rho: \rho\cdot\rho=0\}$ are
uniformly bounded in $\Psi^0(K)$ and $\Psi
^{-1}(K)$, respectively, for $K\subset\subset\R^n$ and hence
interpolation implies
$$\Bigl\{ |\rho|^{1-t} G_\rho:\rho\cdot\rho=0\Bigr\}\subset
\Psi^{-t}(K)\hbox{ is
bounded },\forall t\in [0,1].\tag2.9$$
\bigskip
We can now state a local analogue for the finite-regularity conormal
spaces of
the result of [32] concerning solvability of inhomogeneous equations
involving
$\Dr+q(x)$ in weighted
$L^2$ spaces.
\bigskip
\proclaim{Proposition 2.6} If $q\in I^{\mu}(H)$ with
$\frac23-k\le\mu<1-k$,
$l>3,l\notin\Bbb N$  and $0\le\sigma\le 1$, then for
$\dv_0$ as in (2.8), the inhomogeneous equation
$$(\Dr+q)w=g\in \cl I^{1-k-\sigma,\dv_0}(H)\tag2.10$$
has, for $|\rho|$ large, a unique solution $w\in \cl I^{-k,\dv_0}(H)$,
with $\|w\|\le \frac{C}{|\rho|^\sigma}\|g\|$.\endproclaim
\bigskip
\fp{\bf Proof.} Applying $G_\rho$ to both sides of (2.10),  using (2.9)
for
$t=1-\sigma$ and
Prop. 2.3, we are reduced to showing that
$$(I+G_\rho M_q)w=G_{\rho}g\in\cl  I^{-k,\dv_0}(H)$$
has a unique solution for $|\rho|$ sufficiently large, with $\|w\|\le
C\|G_\rho g\|$.
By Prop. 2.5,
$M_q: \cl I^{-k,\dv_0}(H)\rta \cl I^{\mu+\epsilon,\dv_0}(H)$ for any
$0<\epsilon<-\mu$. Note that
$t=\mu+\epsilon-(-k)<1$, so we can use Prop. 2.5 and (2.9)
with this
value of $t$ to obtain
\medskip
$$\cl I^{-k,\dv_0}(H)\buildrel{M_q}\over{\rta}
\cl I^{\mu+\epsilon,\dv_0}(H)
\buildrel{G_{\rho}}\over{\rta}
\cl I^{-k,\dv_0}(H)$$
with norm
$\le C(q)|\rho|^{t-1+\epsilon'}\rta 0$ as $|\rho|\rta\infty$. Hence, for
$|\rho|$ sufficiently
large, $||G_\rho M_q||<\frac12$ and $I+G_\rho M_q$ is invertible on
$\cl I^{-k,\dv_0}(H)$.\qed
\bigskip
We may now complete the proof of Thm. 1 as described in \S1. Construct two
solutions $\psi_j$, $j=1,2$, to (1.10) for potentials
$q_j\in I^{\mu_j}(H_j)$. Note that $-q_j$, the right hand side of (1.10),
belongs to $I^{\mu_j}(H)\hra\cl  I^{1-k_j-\sigma_j,\dv}(H)$, with
$\sigma_j>0$ since
$\mu_j<1-k_j$. Thus, we may apply Prop.
2.6 with $g=-q_j, j=1,2$, and then form as above the corresponding
solutions,
$v_1(x,\rho_1)=e^{\rho_1\cdot x}(1+\psi_1(x,\rho_1))$ of
$(\D+q_1)v_1=0$ and
$w_2(x,\rho_2)=e^{\rho_2\cdot x}(1+\psi_2(x,\rho_2))$ of
$(\D+q_2)w_2=0$, with $\|\psi_j\|_{\cl I^{-k_j,\dv}}\le
C|\rho|^{-\sigma_j}$. The solutions
$v_1$ and
$w_2$ belong to
$\xpr(\O)$, with
$p,r$ as in (1.3). In fact, each
 is in $L^p(\O)$,  and  in
$W^{2,r}(\O)$  away from $H$ by Prop. 2.2(i) and (iii), resp., since
$r<\frac1{2(1-\nu)}$. Furthermore, since  $q_1\in
L^{\frac{k_1}{\nu_1}-\epsilon},\forall\epsilon>0$, we have  $\D
v_1=-q_1v_1\in
\Bigl(L^{\frac{k_1}{\nu_1}-\epsilon}\Bigr)\times
(L^s\Bigr),\forall\epsilon>0,\forall s<\infty$. Since
$p>\frac{k_1}{k_1-\nu_1}\implies p'<\frac{k_1}{\nu_1}$, we thus
have
$\D v_1\in L^{p'}(\O)$, and similarly for $w_2$. These solutions are
constructed for all large $|\rho|$. Since $n\ge 3$, for any
$\xi\in\R^n\back
0$ and $\lambda\ge c|\xi|$, one can find
$\rho_1,\rho_2\in\{\rho\cdot\rho=0\}$
with
$|\rho_1|\simeq|\rho_2|\ge\lambda$ and $\rho_1+\rho_2=-i\xi$.

By the
assumption that $\CD_{q_2}=\CD_{q_1}$, there exists a $v_2\in\xpr(\O)$
such that
(1.7) holds. Applying (1.8) and (1.9), it suffices to show that
$$\int_\O e^{-i\xi\cdot x}
(q_1-q_2)(\psi_1+\psi_2+\psi_1\psi_2)dx\rta 0\hbox{ as }\l\rta
\infty.\tag2.11$$
Since each $q_j\in L^{\frac{k_j}{\nu_j}-\epsilon}(\O),\forall\epsilon>0$,
and, by (1.3), $p>(\frac{k_j}{\nu_j})', j=1,2$, we have
$$\int_{\O} |q_1-q_2|\cdot |\psi_j|
dx\le\frac{c}{\lambda^{\epsilon_j}},\quad\epsilon_j>0$$
as $\l\rta \infty$ by
H\"older's inequality. In fact,
$p>2(\frac{k_j}{\nu_j})'$ by (1.3), so we  have
$$\int_{\O} |q_1-q_2|\cdot |\psi_1|\cdot |\psi_2|
dx\le\frac{c}{\lambda^{\epsilon_{12}}},\quad \epsilon_{12}>0,$$
as well. Thus,
 we have shown that $\CD_{q_1}=\CD_{q_2}$ implies that
$\widehat{q_1-q_2}(\xi)=0,\forall\xi$, and hence $q_1=q_2$, concluding the
proof  of
Thm.1.

\head
{\bf 3. Reconstruction of the potential}
\endhead

In this section we prove that the potential $q$ can be obtained
constructively from the Cauchy data of $\D+q$.
We follow here the general method of [25]; see also
[26] and [19].
However, there are
some additional difficulties in our case because we deal with the set of
Cauchy data instead of the Dirichlet-to-Neumann map. Moreover, we work
with more complicated function spaces due to the singularities of the
potential. We will show:

\bigskip
\proclaim{Theorem 2} Suppose that $\O$ is Lipschitz and
$H\subset\subset\O$ is a submanifold of codimension $k$.
Suppose further that $q\in I^{\mu}(H)$
is a real potential  with $\nu_0(k)-k<\mu<1-k$ and
$\supp(q)\subset\subset\O$. Let $p,r$ satisfy $2\le
r<r_0(k,\nu)$ and
$\max(p_0(k,\nu),-\frac k\mu)<p<\infty$.
Then $q$ can be reconstructed on
$\O$ from the Cauchy data $\CD_q$ of $\D+q$ on $\xpr(\O)$.
\endproclaim

To start the discussion of reconstruction, we first show how to obtain a
global
analogue of Prop. 2.6. For $s\in\R$ and $-1<\delta<0$, let
$W^{s,2}_\delta(\R^n)$ be the weighted Sobolev space denoted by
$H^s_\delta$ in
[32]. In [32], it is shown that for $0\le t\le 1$,
$$||G_\rho f ||_{W^{s+t,2}_\delta}\le \frac{C}{|\rho|^{1-t}} || f
||_{W^{s,2}_{\delta+1}}.\tag3.1$$

Now, let $\supp(q)\subset \Omega'\subset\subset \Omega$ and let
$\chi_0+\chi_\infty\equiv 1$ be a partition of unity subordinate to the
open
cover $\Omega\cup(\Omega')^c=\R^n$. Fix $\frac23-k<\mu<1-k$, $l>3$ and let
$\dv_0$ be as in (2.8). For $m\le 1-k$ and $s\le -2$, define
$$||f||_{Y^{m,s}_\delta(\R^n)} = ||\chi_0\cdot f ||_{\cl I^{m,\dv_0}(H)}
+ ||\chi_\infty\cdot f||_{W^{s,2}_\delta(\R^n)}.\tag3.2$$
For $m\le -k$, elements of $I^m(H)$ are in $L^p_{loc}(\R^n),\forall
p<\infty$,
and hence $I^m_{comp}(H)\hra W^{s,2}_\delta(\R^n),\forall m\le 1-k, s\le
-1$.
Combining (3.1) with (2.9) and Prop. 2.3, we have, for $0\le t\le 1$,
$$|| G_\rho f ||_{Y^{m-t,s+t}_\delta(\R^n)}\le\frac{C}{|\rho|^{1-t}}||f||
_{Y^{m,s}_{\delta+1}(\R^n)},\quad m\le 1-k-t, s\le -1-t,
-1<\delta<0.\tag3.3$$
Finally, since $\supp(q)\subset\subset\Omega$, it follows from (2.7) that
$$M_q: Y^{-k,s}_\delta\rightarrow Y^{\mu+\epsilon,s}_\delta,\quad\forall
\epsilon>0.\tag3.4$$
Arguing as in the proof of Prop. 2.6,
but replacing the local finite-regularity spaces with the global spaces
$Y^{m,s}_\delta(\R^n)$, we obtain the following result:
\bigskip
\proclaim{Proposition 3.1} If $q\in I^\mu(H)$ with $\frac23-k<\mu<1-k,
s\le -2,
-1<\delta<0$ and $0\le\sigma\le 1$, then the inhomogeneous equation
$$(\Dr+q)w=g\in Y^{1-k-\sigma,s}_{\delta+1}(\R^n)\tag3.5$$
has, for $|\rho|$ large, a unique solution $w\in
Y^{-k,s+1-\sigma}_\delta(\R^n)$ with $||w||\le\frac{C}{|\rho|^\sigma}
||g||$.
\endproclaim
\bigskip

Next we will construct the boundary values of the exponentially
growing solutions on $\partial \Omega$.  For this purpose
we use the Green's function $G_\rho^q(x,y)$ defined by
$$
(\Delta+q)G_\rho^q(\cdotp,y)=\delta_y \quad\hbox{in}\quad \R^n,\ \ 
 e_\rho(\cdotp) G_\rho^q(\cdotp,y)\in Y^{-k,s}_{\delta}(\R^n), \tag3.6
$$
where $y\in \R^n\setminus \overline\Omega$,
$e_\rho(x)=\exp(-\rho\cdotp x)$ and $s<-n/2$. When $|\rho|$
is large enough, the equation (3.6) has a unique
solution by Prop. 3.1.
Next we consider the case when
$\rho$ is fixed and sufficiently large.

 As $\supp(q)\subset\subset \O$, we see that $\p\O$ has a neighborhood $V$
such that in $V\times V$
Green's function $G_\rho^q(x,y)$ has the
same singularities as the  Green's function $G^0_0$ (for the zero
potential and $\rho=0$),
that is, $G_\rho^q(x,y)-G_0^0(x,y)\in C^\infty(V\times V)$.

Using the Green's function (3.6) we
define the corresponding single and double layer potentials
$$
S_{q}\phi(y)=\int_{\p\O}G_\rho^q(x,y)\phi(x)dS_x,\quad
K_{q}\phi(y)=\int_{\p\O} (\frac \p{\p n(x)}
G_\rho^q(x,y))\phi(x)dS_x,\quad y\not\in \p\O
$$
which define continuous operators $S_{q}:W^{1-\frac 1r,r}(\p\O)\to
X|_\O \oplus
W^{2,r}_{loc}(\R^n\setminus\O)$
and $K_{q}:W^{2-\frac 1r,r}(\p\O)\to
X|_\O \oplus W^{2,r}_{loc}(\R^n\setminus\O)$. Here
$X\subset {\Cal D}'(\R^n)$ is the space with the norm
$||\chi_0\cdot f ||_{\cl I^{m,\dv_0}(H)} +
||\chi_\infty\cdot f||_{W^{2,r}_\delta(\R^n)}.$
These layer potentials can be considered
as operators on the boundary $\p\O$, defined in principal value
sense. Since $\p\O$ is Lipschitz, it follows from the results of
[7] that these operators are continuous,
$S^{\p\O}_{q}:W^{1-\frac 1r,r}(\p\O)\to
W^{2-\frac 1r,r}(\p\O)$
and
$K^{\p\O} _{q}:W^{2-\frac1r,r}(\p\O)\to
W^{2-\frac 1r,r}(\p\O)$.
 Similarly, on $\p\O$ we
define normal derivatives of the layer potentials,
$$
\p_n S^{\p\O}_{q}\phi(y)=\hbox{p.v.}\int_{\p\O} (\frac \p{\p n(y)}
G_\rho^q(x,y)\phi(x)dS_x,
$$
$$
\p_n K_{q}^{\p\O}\phi(y)=\hbox{p.v.}\int_{\p\O} (\frac \p{\p n(y)}
 \frac \p{\p n(x)}
G_\rho^q(x,y))\phi(x)dS_x
$$
which are continuous operators
$\p_n S^{\p\O}_{q}:W^{1-\frac 1r,r}(\p\O)\to
W^{1-\frac 1r,r}(\p\O)$
and
$\p_n K^{\p\O}_{q}:W^{2-\frac 1r,r}(\p\O)\to
W^{1-\frac 1r,r}(\p\O)$.

Next we consider the Calder\'on projector [5]. We
start with the operator
$$
A_q(
\phi,\psi)=
\left(- S^{\p\O}_q\phi+(-\frac 12+K^{\p\O}_q)\psi,
-(\frac 12+\p_n S_q^{\p\O})\phi+\p_n K_q^{\p\O}\psi\right).
$$
\medskip
\proclaim{Proposition 3.2} Let $Z=W^{2-\frac 1r,r}(\p\O)\times W^{1-\frac
1r,r}(\p\O)$. Then the operator
$$
A_q:Z/
\hbox{Ker}(A_q)\to Z
$$
is semi-Fredholm. Moreover, $-A_q:Z\to Z$ is a projection
operator with range
$\CD_q$ and kernel independent of $q$. In particular,
$\CD_q$ is a closed
subspace of $Z$.
\endproclaim
\medskip
\fp Proof.
First  we show that kernel of $A_q$ does not depend on $q$.
Assume that $(\phi,\psi)\in Z$. We consider the
function
$$
u_{\phi,\psi}(y)=-S_q(\phi)+K_q(\psi),\quad y\in \R^n\setminus \p \O
$$
and the trace operators
$$
T_+:W^{2,r}(\R^n\setminus \O)\to W^{2-\frac 1r,r}(\p\O)\times
W^{1-\frac 1r,r}(\p\O),
$$
$$
T_-:W^{2,r}(\O)\to W^{2-\frac 1r,r}(\p\O)\times
W^{1-\frac 1r,r}(\p\O)
$$
defined by $T_\pm u=(u|_{\p\O},\p_n u|_{\p\O})$.
As the Green's functions $G^q_\rho(x,y)$ have
the same singularities near $\p \O\times \p\O$
as the standard Green's function of $\R^n$,
we can use the standard
jump relations for layer potentials (see e.g. [8]). We conclude that
$$
T_- u_{\phi,\psi}=A_q(\phi,\psi),\quad
T_+ u_{\phi,\psi}=(\phi,\psi)+A_q(\phi,\psi).
$$
Thus we get that $u=u_{\phi,\psi}\in  (e_\rho)^{-1}Y_\delta^{-k,s}$ and
it is the unique solution of
$$
(\Delta+q)u=g_{\phi,\psi}=\psi\delta_{\p \O}+
\nabla\cdotp (n\phi
\delta_{\p \O}) \quad\hbox{in}\quad \R^n, \tag3.7
$$
satisfying
$ e_\rho(\cdotp) u\in Y^{-k,s}_{\delta}(\R^n).$

Now, if $(\phi,\psi)\in\hbox{Ker}(A_q)$ we have
that
$$
(\phi,\psi)=(\phi,\psi)+T_- u_{\phi,\psi}=T_+ u_{\phi,\psi}.
$$
Thus, $v=u_{\phi,\psi}$ is
the solution of the scattering problem
$$
\Delta v=0 \quad\hbox{in}\quad \R^n\setminus\O,\ \
T_+ v=(\phi,\psi),\quad
e_\rho(\cdotp) v\in W^{s,2}_{\delta}(\R^n\setminus \O).
\tag3.8
$$
On other hand, assume that (3.8) has a solution, and
let $v_0$ be the zero-continuation of $v$,
that is $v_0|_{\R^n\setminus \O}=v,$ $v|_\O=0$. Then
we conclude that $v_0$ is a solution
of the problem (3.7), and as this solution is unique,
we see that $v_0=u_{\phi,\psi}$. This shows that
$(\phi,\psi)\in\hbox{Ker}(A_q)$ if and only if
the problem (3.8) has a solution. This is obviously independent of $q$
and thus we see that
$$
T_-+(e_\rho+u_{\phi,\psi})=(\phi_\rho+\phi,\psi_\rho+\psi)+
A_q(\phi,\psi)\in\hbox{Ran}(A_q).\tag3.9
$$
Applying the projection $I+A_q$ to both sides of (3.9) and using
$A_q(\phi,\psi)=0$, we see that
$$
0=(I+A_q)(\phi_\rho+\phi,\psi_\rho+\psi)=(I+A_q)(\phi_\rho,\psi_\rho)+(\phi,\psi).
$$
As $A_q$ and $(\phi_\rho,\psi_\rho)$ are known, we can thus determine
$(\phi,\psi)$ and the Cauchy data of $v(x)$ on $\p\O$.

Next we consider range of $A_q$. A standard application
of Green's formula (see, e.g., [8, Th. 3.1]) shows that if $v\in X^{p,r}$
satisfies
$$
(\Delta+q)v=0 \quad\hbox{in}\quad \O,
$$
and $(\phi,\psi)=T_-v$, then $v=-u_{\phi,\psi}.$ (Observe
the negative sign which is due to the fact that we use
exterior normal vector $n$.)
Also, by for $(\phi,\psi)\in Z$ we have
$\chi_0u_{\phi,\psi}\in
\cl I^{-k,\dv_0}(H)\subset L^{t_1}(\O)$ for any $t_1<\infty$
by Prop. 2.2. As $q\in I^{\mu}(H)\subset L^{t_2}(\O)$
for $1<t_2<\frac k{k+\mu}$ we have
$\chi_0\Delta u_{\phi,\psi}=-\chi_0 qu_{\phi,\psi}\in L^{p'}(\O)$
for $\frac 1p <1-\frac 1{t_2}$, i.e. $p>-\frac k\mu$. Hence
$u_{\phi,\psi}\in X^{p,r}$. Thus the set of all solutions
of the Schr\"odinger equation in $X^{p,r}$ equals to the set of
solutions $u_{\phi,\psi},$ $(\phi,\psi)\in Z$.
As $T_-u_{\phi,\psi}=A_q(\phi,\psi)$,
we obtain that the range of $A_q$ equals to $\CD_q$.

Now, when the potential
is equal to zero, the
Dirichlet-to-Neumann operator $\Lambda_0:u|_{\p\O}=\p_n u|_{\p\O}$
is well defined, $\Lambda_0:W^{2-\frac 1r,r}(\p\O)\to
W^{1-\frac 1r,r}(\p\O)$. The Cauchy data $\CD_0$ is the graph of the
operator $\Lambda_0$ and  is thus closed.
Thus we see that the range of $A_0$ is a closed subspace,
and therefore the operator
$
A_0:Z/\hbox{Ker}(A_0)\to Z
$
has zero kernel and closed range. Thus it is
a semi-Fredholm operator. Now, consider the operator
$A_q-A_0$. Using (3.9) we know that the operator
$$
A_q-A_0:Z/\hbox{Ker}(A_0)\to Z
$$
is well defined and compact. As compact perturbations
of semi-Fredholm operators are also  semi-Fredholm,
we conclude that $A_q$ is also semi-Fredholm.

It remains to show that $-A_q$ is a projection.
This can be seen similarly to the  smooth case. Indeed,
if $(\phi,\psi)\in \hbox{Ran}(A_q)$,
$(\phi,\psi)=A_q(\tilde \phi,\tilde \psi)$ we see that
the solution $u_{\tilde \phi,\tilde \psi}$ has the trace
$T_-u_{\tilde \phi,\tilde \psi}=(\phi,\psi).$ Hence
Green's formula gives
$$
u_{\tilde \phi,\tilde \psi}=-(-S_q\phi+K_q\psi)\quad\hbox{in}\quad \O.
$$
Taking trace $T_-$ from both sides we
obtain that $(\phi,\psi)=-A_q(\phi,\psi)$, i.e. $(-A_q)^2=-A_q$.
Thus, Prop. 3.2 is proven.\qed
\bigskip
Now we can construct the boundary values of the
exponentially growing solutions from the Cauchy data. As we are given
$\CD_q=\hbox{Ran}(A_q)$, and we know
$\hbox{Ker}(A_q)=\hbox{Ker}(A_0)$, we can construct
the projection $-A_q$. Next, let
$(\phi_\rho,\psi_\rho)=
T_+e_\rho$ be the boundary values of the incoming plane
wave. Consider the solution
$v(x)= e^{\rho\cdot x}(1+\psi(x,\rho))=e^{\rho\cdot x}+u_0$ and let
$(\phi,\psi)=T_+u_0$. Then $(\phi,\psi)\in \hbox{Ker}(A_q)$
and $u_0=u_{\phi,\psi}$ in $\R^n\setminus \O$. Moreover, as $v$ is
solution
of Schr\"odinger equation inside $\O$, we have
$$
T_+(e_\rho+u_{\phi,\psi})=(\phi_\rho+\phi,\psi_\rho+\psi)+
A_q(\phi,\psi)\in \hbox{Ran}(A_q).\tag3.9
$$
Applying with projection $I+A_q$ to (3.9) and using $A_q(\phi,\psi)=0$,
we see that
$$
0=(1+A_q)(\phi_0+\phi,\psi_0+\psi)=
(1+A_q)(\phi_0,\psi_0)+(\phi,\psi).
$$
As $A_q$ and $(\phi_0,\psi_0)$ are known, we find $(\phi,\psi)$
and the Cauchy data of $v(x)$ on $\p \O$.

So far, we have constructed the Cauchy data of the
solutions $v_{\rho_1}(x)= e^{\rho_1\cdot x}(1+\psi(x,\rho_1))$
for all sufficiently large $\rho_1$. Thus if we consider
complex frequencies $\rho_1$ and $\rho_2$
satisfying
$\rho_1+\rho_2=-i\xi$, with $\xi\in\R^n\back 0$, an application of Green's
formula yields
$$
\eqalign{{\hat q}(\xi)&=\lim_{|\rho_1|\to \infty}
 \int_\O q(x) e^{\rho_1\cdot x}
\Bigl(1+\psi_1(x,\rho_1)\Bigr)\cdotp e^{\rho_2\cdot x}\,dx
\cr
=&\lim_{|\rho_1|\to \infty}
 \int_{\p \O}\Bigl( v_{\rho_1}\cdotp \p_ne^{\rho_2\cdot x}
- \p_n v_{\rho_1}\cdotp e^{\rho_2\cdot x}
\Bigl)\,dx
\cr}
$$
This proves Theorem 2.\qed

\head
{\bf 4. Non-uniqueness for highly singular potentials}
\endhead

We next discuss how very strong singularities of the
potential can
cause non-uniqueness in a closely related inverse problem. Due to the
strength of the singularities,  the Schr\"odinger equation has to be
interpreted in a weak sense. Let us consider the boundary value problem
$$
(\Delta+q+E)u=0\quad \hbox{in}\quad \Omega,\ \ \ u|_{\p\O}=f
\tag4.1
$$
with the potential $q$ having the form
$$
q(x)=-\dist(x,H)^{\mu}c_0(x),\tag 4.2
$$
near $H$, where $\dist$ is the Euclidean distance, $H$ is a closed
hypersurface bounding a region $\O_0\subset\subset\O$,
$\mu<-2$, and $c_0(x)$ is a smooth function,
satisfying
$
c_0(x)>C_0>0
$
in some neighborhood $V$ of $H$.

 Elements of  $ I^{-1-\mu}(H)$ satisfy the pointwise estimate
$|q(x)|\le C \dist (x,H)^\mu$, but a $q$ satisfying (4.2) is not
even locally integrable and thus need not define a distribution. Hence,
the solutions of (4.1) cannot be formulated in the usual sense of
distributions. Instead,  we define the solution of (4.1) (if it
exists) to be
the solution of the following convex minimization problem:
Find $u$ such that
$$
G(u)=\inf G(v)\tag4.3
$$
where $G=G_{q+E}:\{v\in H^1(\O):\ v|_{\p\O}=f\}\to \R\cup\{\infty\}$
is the convex functional
$$
G_{q+E}(v)=\int_\O (|\nabla v(x)|^2-(q(x)+E)
|v(x)|^2)dx
$$
Here, since the  function $-q(x)$ is bounded from below, we define
$G(u)=\infty$ when $q|v|^2$ is not in $L^1(\O)$.

\bigskip
\proclaim{Proposition 4.1}
The Cauchy data
$$
\CD_{q+E}=
\Big\{(u|_{\po},\frac{\partial u}{\partial n}|_{\po}): u\in
H^1(\O),\ u\hbox{ is a minimizer of }G_{q+E}\}.
$$
does not depend on $q$ in $\O_0$. In particular, if
the solution of (4.3) is unique, $u$ vanishes identically in $\O_0$.
\endproclaim
\bigskip
{\bf Remark.} We note that potentials having singularities similar to
(4.2) as above  has been used
to produce counterexamples to strong unique continuation, e.g.
potentials $q(x)=c/\vert x\vert \sp {2+\varepsilon}$ in
[11].  Recently,
counterexamples  have been found for weak unique continuation for
$L^1$-potentials [18], but here we need to construct potentials for which
{\it all} solutions vanish inside $H$.
 Finally,
we wish to emphasize that since the solutions of (4.1) considered here
are not defined in the usual sense of distributions,
but rather as solutions of a convex minimization problem, the solutions
 we construct do not give new counterexamples for the unique
continuation problem.

\smallskip

\fp{\bf Proof.} Obviously we can assume that $q(x)\leq 0$ everywhere.
We start first with the case where $E=0$ and $f\in C^\infty(\p \O)$.

As the potential $q$ is not in the Kato class ([6,p.62]),
 consider instead  a decreasing sequence of  smooth functions
 $q_n\in C^\infty(\O)$, $q_{n+1}(x)\leq q_n(x)$,
for which $q_n(x)=q(x)$ when
$d(x,H)>\frac 1n$ and
in some neighborhood $V$ of $H$
$$
q_n(x)\leq \max(-c_1 n^{-\mu},q(x))\tag4.4
$$
where $0<c_1<C_0$.
Let $G_n$ be the functionals defined as
$G$ with $q$ replaced with $q_n$. The functionals $G_n$ have
unique minimizers $u_n$ which satisfy in classical sense
$$
(\Delta+q_n)u_n=0\quad \hbox{in}\quad \Omega, \ \ u_n|_{\p\O}=f.\tag4.5
$$
 Now, let $f\in C^\infty(\p\O)$ be fixed.
Let $F\in H^1(\O)$ be a function for which $F|_{\p\O}=f$ and $F=0$ in
some neighborhood of $H$. By definition of the potentials $q_n$,
for sufficiently large $n_0$ we
have $G(F)=G_n(F)=G_{n_0}(F)$ for $n\geq n_0$. Thus
for the minimizers $u_n$ of $G_n$ we have
$
G_n(u_n)\leq C=G_{n_0}(F).
$
Next, by choosing a subsequence, we can
assume that the sequences $\int|\nabla u_n(x)|^2dx$ and $ \int(-q_n(x))
|u_n(x)|^2dx$
are decreasing when $n\to\infty$. Next,
let us denote by $C_1,C_2\leq C$ the constants
$$
C_1=\lim_{n\to\infty}
\int_\O |\nabla u_n(x)|^2dx,\quad
C_2=\lim_{n\to\infty} \int_\O (-q_n(x)) |u_n(x)|^2dx=C_2.
$$

Now, we see that $u_n$ are uniformly bounded in $H^1(\O)$
and thus by choosing a subsequence we can assume
that there is $\tilde u\in H^1(\O)$ such that
 $u_n\to \tilde u$ weakly in $H^1(\O)$. Moreover,
$$
\int_\O |\nabla \tilde u(x)|^2dx\leq C_1 \tag4.6
$$
As compact operators map weakly
converging sequences to strongly converging ones, we have
the norm-convergence
$u_n\to \tilde u$ in $L^2(\O)$. Thus for $n_0>0$ we have
$$
-\int_\O q_{n_0}(x) |\tilde u(x)|^2dx
=\lim_{n\to\infty} \int_\O (-q_{n_0}(x)) |u_n(x)|^2dx\leq
$$
$$
\leq
\lim_{n\to\infty} \int_\O (-q_{n}(x)) |u_n(x)|^2dx\leq
C_2.
$$
Since this is valid for any $n_0$ we have by monotone convergence
theorem
$$
-\int_\O q(x) |\tilde u(x)|^2dx \leq
C_2. \tag4.7
$$
As $G\leq G_n$,
$$
\inf G(u)\leq \lim_{n\to\infty} \min G_n(u)=C_1+C_2.
$$
Now, by (4.6) and (4.7) we have
$
G(\tilde u)\leq C_1+C_2
$
and thus $\tilde u$ is a minimizer of $G$. Now,
as $G=G_n+(G-G_n)$ where $G_n$ is a strictly convex functional
and $G-G_n$ is a convex functional, $G$ is strictly convex.
Thus the minimizer is unique. Hence we see that, for
every $f$, the solution $\tilde u$ of the minimization problem
(4.3) exists, is unique, and is given as the $L^2$-limit
of the functions $u_n$. We note that the above
analysis was based to the fact that the minimization problems
for the $G_n$ epi-converge to the minimization problem for $G$
[30].
\bigskip
Recalling that $\O_0$ is the region bounded by $H$,
consider functions $u_n$ restricted to $\O_0$.
Let  $t\mapsto B_t$ be the Brownian motion in $\R^n$
starting from $x$ at time $t=0$, i.e., $B_0=x$.
As the $q_n$ are strictly negative smooth functions, they
are in the Kato class and the pair $(\O,q_n)$ is gaugeable
(see [6], sect. 4.3 and Th. 4.19).
By [6], Th. 4.7,
the solution $u_n$ can be represented by the Feynman-Kac formula
$$
u_n(x)=E(\exp\left(\int_0^\tau q_n(B_t)dt\right)f(B_\tau))
$$
where $\tau=\tau_{\p\O}$ is the first
time when the process hits the boundary, i.e., $B_t\in \p \O$.
Here, we assume $B_t$ is a version of Brownian motion for
which all realizations are continuous curves (see [21] or [6], Th.
1.4).  If $x\in \O_0$, the realizations of Brownian motion have
to hit $H$ prior to hitting  $\p \O$. Denote
the first hitting time for $H$ by $\tau_H$; thus the first
hitting point is $B_{\tau_H}$, and $\tau_H<\tau_{\p\O}$. (The stopping
time
$\tau_H$ is measurable function in the probability space,
 see [6,Prop. 1.15]).

Let us now denote by
$p(\rho,\eta)$ the probability that the Brownian motion
sent from origin at time $t=0$ leaves
the origin centered ball with radius $\rho$ before time $\eta$.
Because of the scale-invariance
of Brownian motion, $p(s\rho,s^2\eta)=p(\rho,\eta)$ for $s>0$.
(Indeed, let us consider reparametrized Brownian motion
$\tilde B_t=sB_{s^{-2}t}$. As the probability densities
of $(\tilde B_{t_1},\tilde B_{t_2},\dots,\tilde B_{t_m})$
coincide to those of  $(B_{t_1},B_{t_2},\dots,B_{t_m})$
we see that we see that $\tilde B_t$ is Brownian motion, too.)

Let $A_{\rho,\eta}=\{|B_t-B_{\tau_H}|<\rho\hbox{ for }\tau_H\leq
t<\tau_H+\eta\}$.
This set is measurable in the probability space and
the probability of
$A_{\rho,\eta}$ is
$P(A_{\rho,\eta})=1-p(\rho,\eta)$.

Let $m>1$  and
$\eta=\eta(m)$ be such that
$p(1,\eta)\geq
\frac {m-1}m$. Now, $q$ is non-positive and  by (4.4)
$q_n(x)<\max(-c_1n^{-\mu},q(x))$ in some neighborhood $V$
of $H$. When $s$ is so small that the $s$-neighborhood of $H$
is in $V$, we have by (4.2) that
$$\eqalign{
&|E(\exp\left(\int_0^\tau q(B_t)dt\right)f(B_\tau))|\leq
E(\exp\left(\int_{\tau_H}^{\tau_H+s^2\eta}
q(B_t)dt\right)||f||_{L^\infty})\cr
\leq&
(1-P(A_{s\rho,s^2\eta(m)}))||f||_{\infty}
+
P(A_{s\rho,s^2\eta(m)})
\exp\bigg(-s^2\eta(m)
\min(C_0s^\mu,c_1n^{-\mu})\bigg)||f||_{\infty}.}
$$
Thus, choosing $s=n^{2/\mu}$ we see that for sufficiently large $n$
$$
||u_n||_{L^\infty(\O_0)}\leq (\frac 1m+\frac {m-1}m
\exp\left(-\eta(m)c_1n^{2(2+\mu)/\mu}\right)||f||_{L^\infty}. \tag4.8
$$
As $u_n\to \tilde u$ in norm in $L^2(\O_0)$,
$
||\tilde{u}||_{L^2(\O_0)}\leq \frac 1m
||f||_{L^\infty}\,\hbox{vol}(\O_0)^{1/2}
$
for any $m$. Thus we see that
$\tilde u=0$ in $\O_0$.

Next we consider the case when $E\in \R$ and $f\in H^{1/2}(\p\O)$.
First, let
$H_r=\{x\in \O:\ \dist(x,H)<r\}$ and let $r$ be so small
that $q(x)+E<0$ for $x\in H_r$. If $u$ is the solution of (4.1)
in $\O$, then its restriction $\tilde u=u|_{H_r}$
is the solution of boundary value problem
$$
(\Delta+q+E)\tilde u=0\quad \hbox{in}\quad H_r,\ \ \tilde u|_{\p H_r}=\tilde f, \tag4.9
$$
where $\tilde f=u|_{\p H_r}\in C^\infty (\p H_r),$
that is, $\tilde u$ is the solution of minimization problem
(4.3) in domain $H_r$.
Let $q_n$ approximate $q$ in $H_r$ as above and $\tilde u_n$ be
the corresponding solutions of problem (4.9) with $q$ replaced with $q_n$.
As above, we see that problem (4.9)
is uniquely solvable,
 $\tilde u_n\to \tilde u$ in $L^2(H_r)$, and that
 $\tilde u_n(x)$ can
be represented by the Feynman-Kac formula. Let $x\in \O_0\cap H_r$,
$\tilde \tau$ be the first time when the Brownian motion
sent from $x$ at $t=0$ hits  $\p H_r$,
and $\tilde A=\{B_{t}\in \O_0\cap H_r\hbox{ for }0\leq t<\tilde \tau\}$.
Let us denote $\tilde f=\tilde f_++\tilde f_-$,
where $\tilde f_+$ vanishes on $\O_0\cap \p H_r$ and
 $\tilde f_-$ vanishes on $(\O\setminus \O_0)\cap \p H_r$.
Then we see that
$$
\tilde u_n(x)=P(\tilde A)E(\exp\left(\int_0^{\tilde \tau}
q_n(B_t)dt\right)
\tilde f_-(B_{\tilde \tau})|\tilde A)+
$$
$$
+(1-P(\tilde A))E(\exp\left(\int_0^{\tilde \tau}
q_n(B_t)dt\right)\tilde f(B_{\tilde \tau})|\tilde A^c)
\tag4.10
$$
where $E(\cdot|\tilde A)$ is conditional expectation with condition
$\tilde A$ and $\tilde A^c$ denotes the complement of $\tilde A$.
Note that in the case of $\tilde A^c$, the process $B_t$ hits
 $H$ at least once.
Analyzing how long Brownian motion is near $H$
as above, we see that when $n\to \infty$,
the second term on the right hand side of (4.10)
goes to zero. Thus $\tilde u(x)$, for $x\in H_r\cap \O_0$, depends
only on $q$ in $H_r\cap \O_0$ and $f_-$.
Similarly, we
see that  $\tilde u(x)$, $x\in H_r\setminus\overline \O_0$, depends
only on $q$ in $H_r\setminus\overline \O_0$ and $f_+$.
Moreover, analogously to (4.8)
we see that
$$
\lim_{n\to\infty}||\tilde u_n|_H||_{L^\infty(H)}=0.
$$
Choosing a subsequence, we can assume that
 $\tilde u_n\to \tilde u$ weakly in $H^1(H_r)$
and thus in norm in  $H^{3/4}(H_r)$. Hence, by taking the
trace $H^{3/4}(H_r)\to L^2(H)$ we see that
$\tilde u|_H=0$.

In conclusion, for the boundary value problem (4.9) there are
well defined maps
$$
T_+:f_+\mapsto \tilde u|_{H_r\setminus \overline \O_0}\in
\{v\in H^1(H_r\setminus \overline \O_0):\ v|_H=0\},
$$
$$
T_-:f_-\mapsto \tilde u|_{H_r \cap \O_0}
\in \{v\in H^1(H_r \cap \O_0):\ v|_H=0\}
$$
where $T_+$ depends only on $q$ in $H_r\setminus \overline \O_0$
and $T_-$ on  $q$ in $H_r\cap\O_0$.

In particular, on the boundaries
$\p H_r\cap \O_0$ and
$\p H_r\cap (\O\setminus \O_0)$ we have ``independent''
Dirichlet-to-Neumann maps
$$
\Lambda_+:f_+\mapsto \p_n\tilde u|_{\p H_r\setminus \overline \O_0},\quad
\Lambda_-:f_-\mapsto \p_n\tilde u|_{\p H_r \cap \O_0},
$$
where $n$ is the exterior normal of $H_r$.

Next, if $u$ is a solution of boundary value
problem (4.1) we denote
$u_+=u|_{\O \setminus \overline \O_0}$ and
$u_-=u|_{\O_0}$. To motivate the next step, we observe that
$u_+$ and $u_-$ satisfy ``independent'' boundary value
problems in $\O\setminus (H_r\cup \O_0)$
$$
(\Delta+q+E)u_+=0,\ \
u_+|_{\p H_r}=\tilde f,\ \
\p_n u_+|_{\p H_r\setminus \O_0}=\Lambda_+(u_+|_{\p H_r\setminus \O_0})
$$
and in $ \O_0\setminus H_r$ 
$$ 
(\Delta+q+E)u_-=0, \ \ 
\p_n u_-|_{\p H_r\cap \O_0}=\Lambda_-(u_-|_{\p H_r\cap\O_0}).
$$

Now, considering the form of $G$ and the fact that
the solution $u$ of boundary value
problem (4.1) satisfies $u|_H=0$, we see that
$u_+=u|_{\O \setminus \overline \O_0}$
is a minimizer of $G$ in the set
$\{v\in H^1(\O \setminus \overline \O_0):\ v|_{\p \O}=f,\ v|_H=0\}$ and
$u_-=u|_{\O_0}$ is a minimizer of $G$ in the set
$\{v\in H^1(\O_0):\ v|_{\p \O_0}=0\}$.

Conversely, if
$U=v_+$ in $\O \setminus \overline \O_0$ and $U=v_-$ in $\O_0$
where $v_+$ and $v_-$ are any minimizers of $G$
in the sets
$\{v\in H^1(\O \setminus \overline \O_0):\ v|_{\p\O}=f,\ v|_H=0\}$ and
$\{v\in H^1(\O_0):\ v|_{\p \O_0}=0\}$, respectively,
then $U$ is solution of (4.1).

In particular, we see that the Cauchy data of solutions $u$
of (4.1)
on $\p\O$ are independent of $u|_{\O_0}$ and thus of $q$ inside $H$.
This finishes the proof of Prop. 4.1. As a concluding remark we note
that by using the Courant-Hilbert min-max principle, we
see that there always are values of $E$ such that
minimization problem for $v_-$ has non-zero solutions, that is,
there are eigenstates $U$ which have vanishing Cauchy data
on $\p \O$.
\qed
\medskip

Physically, this example has the following interpretation:
In theory it is possible to construct a potential wall $q(x)$
such that
no particles can ``tunnel" through it, using an analogy with quantum
mechanics. Thus
exterior observers can make no conclusions about
the existence of objects or structures
inside this wall. Moreover, inside $H$ the solution can be
in an eigenstate and its Cauchy data vanishes on the boundary of $\O$.
Thus, making another analogy with quantum mechanics, in this nest the
Schr\"odinger cat could live happily forever.

\enddocument